\documentclass{amsart}

\usepackage{amssymb,amsmath,amsthm,amscd}

\newcommand{\N}{\mathbb{N}}
\newcommand{\R}{\mathbb{R}}
\newcommand{\C}{\mathbb{C}}

\newtheorem{theo}{Theorem}[section]%
\newtheorem{cor}[theo]{Corollaire}%
\newtheorem{prop}[theo]{Proposition}%
\newtheorem{rem}[theo]{Remarque}%
\newtheorem{defn}[theo]{Definition}%

\begin{document}

\title{Disques J-holomorphes contenus dans une hypersurface}

\author[E.Mazzilli]%
{Emmanuel Mazzilli}%


\address{E.M.: UFR de Math\'ematiques\newline%
 \indent Universit\'e de Lille 1\newline%
 \indent 59655 Villeneuve d'Ascq\newline%
 \indent France} %
\email{mazzilli@math.univ-lille1.fr}

\begin{abstract}\hskip5pt
We study germs of J-Holomorphic curves
contained in $M$, a real analytic hypersurface of an symplectic
manifold of dimension 4. We show, under topological hypothesis on
$M$, that if $M$ is compact then $M$ is of finite type and so there
is no germs of $J$-holomorphic curves on $M$(with $J$ adapted with
the symplectic form). In $\Bbb{C}^2$ with the standard complex
structure, this is a classical result of Diederich-Fornaess.
\end{abstract}

\maketitle

\section{Introduction et principaux r\'esultats}

Consid\'erons, $(\R^{4},J)$ avec $J$ une structure presque
complexe r\'eelle analytique, $M=\{\rho=0\}$ un germe d'hypersurface
r\'eelle analytique en $0$. Dans [1], nous avons d\'evelopp\'e une
th\'eorie du type r\'egulier dans le cas non int\'egrable; pour
\'enoncer pr\'ecis\'ement certains r\'esultats, il faut rappeler
quelques d\'efinitions et \'enonc\'es de [1] :
\begin{defn}Soit $M=\{\rho=0\}$ un germe
d'hypersurface lisse en $0$; on dit que $0$ est de type fini s'il
existe $k\in \N$ tel que : $\forall \gamma$ disque $J$-holomorphe
r\'egulier $z\rightarrow \gamma(z)$, $0\rightarrow 0$, l'ordre
d'annulation de $\rho o \gamma$ est inf\'erieur \` a $k$. S'il
n'existe pas de $k$ ayant la propri\'et\'e ci-dessus, on dit que $0$
est un point de type infini.\end{defn}
\begin{defn}Si $0$ est de type fini, on appelle type
de $M$ en $0$, le plus petit des $k$ v\'erifiant la d\'efinition 0-1
et l'on note $\delta_M(0)$.\end{defn}
On \'etablit dans [1], la proposition suivante:
\begin{prop}Si $M\subset  \R^{4}$, alors $z\in
M\rightarrow \delta_M(z)$ est une fonction semi-continue sup\'
erieurement.\end{prop}Si $0$ est de type infini, il semble
raisonnable d'esp\'erer trouver un germe de disque $J$-holomorphe,
$\gamma$, contenu dans $M$ ; dans le cas d'une structure
int\'egrable, ceci est vraie, voir par exemple [2], ou bien en
remarquant qu'alors, la vari\'et\'e de Segre de $M$ en $0$ a un
ordre de contact infini avec $M$, (gr\^ace \`a une application de
l'in\'egalit\'e de Lojasciewicz) et donc $\rho(0,z)\subset M$.

Quoiqu'il en soit, l'une et l'autre des m\'ethodes utilisent
fortement l'int\'egrabilit\'e de la structure (par exemple, il n'y a
pas l'equivalent des vari\'et\'es de Segre dans le cadre preque
complexe), et donc ne s'adaptent pas \`a notre situation. Pour
contourner cette difficult\'e, nous sommes amen\'es \` a utiliser
une version ad hoc (voir[4]et [5]) du th\'eor\` eme de
Cartan-Kahler, en remarquant que $0$ de type infini entra\^\i ne la
``formelle int\'egrabilit\'e'' d'un certain syst\` eme
diff\'erentiel (que l'on imagine ais\' ement, voir partie 3). On
montre donc le th\'eor\` eme et le corollaire suivants :
\begin{theo}\label{thm:main}Soit $0\in M\subset \R^{4}$, un germe
d'hypersurface r\'eelle analytique. Alors si $0$ est de type infini,
il existe un disque $\gamma$ contenu dans $M$, avec $\gamma(0)=0$.
\end{theo}

Soit $E$ l'ensemble suivant:
$$E:=\{X\in M/\ \exists \gamma,\ \hbox{disque holomorphe non trivial v\'erifiant :}\
\gamma(0)=X\ \hbox{et}\ \gamma\subset M\}.$$ \begin{cor}
Sous les hypoth\` eses du r\'esultat pr\'ec\'edent, $E$ est
ferm\'e dans $M$.
\end{cor}
Preuve du corollaire 0-5 :

Soit $(a_n)$ une suite de points de $E$ qui converge vers $0\in  M$
et $(\gamma_n)$ la suite des germes associ\'es; si $\gamma_n$ est
singulier en $a_n$, alors il existe une suite de points $(b_k)$ de
$\gamma_n$ qui converge vers $a_n$, telle que $\gamma_n$ r\'egulier
en $b_k$. D'apr\`es la proposition 0-3, $a_n$ est donc un point de
type r\'egulier infini et par cons\'equent le th\'eor\` eme 0-4
entra\^\i ne l'existence de $\beta_n\subset M$ passant par $a_n$ et
r\'egulier. A nouveau, la proposition 0-3 assure que $0$ est un
point de type r\'egulier infini, et le th\'eor\` eme 0-4 qu'il
existe $\beta\subset M$ r\'egulier passant par $0$.

En fait, il n'est point n\'ecessaire de s'assurer que les germes
associ\'es aux points de $E$ sont r\'eguliers; si $M$ est une
hypersurface dans $\R^{4}$, elle ne peut contenir de germes de
disques singuliers (mais pour obtenir ce r\'esultat, il faut une
\'etude plus fine des disques $J$-holomorphes contenus  dans une
hypersurface de $\R^{4}$; voir section 4).

\begin{defn}On dit qu'une hypersurface $M$ dans une
vari\'et\'e presque complexe $(V,J)$ est de type fini,s'il existe
$k$ tel que tout point de $M$ est de type inf\'erieur \` a
$k$.\end{defn}D'apr\` es la proposition 0-3, si $M$ est compacte
et  n'est pas de type fini, alors il existe $p_0\in M$ tel que $p_0$
est de type infini.

Nous avons \'egalement le r\'esultat global suivant :

\begin{theo}\label{thm:main}
 Soit $(V,\omega)$ une vari\'et\'e
symplectique r\'eelle analytique de dimension 4 avec $\omega$ exacte
sur $V$. Alors si $J$ est une structure presque complexe adapt\'ee
\`a $\omega$, toute hypersurface compacte $M$ est soit feuillet\'ee,
soit de type fini.
\end{theo}
\begin{cor} Sous les m\^emes hypoth\` eses que le
th\'eor\` eme pr\'ece\` edent et si de plus $M$ v\'erifie l'une des
deux assertions suivantes : $\pi_{2}(M)\not =0$ ou $\pi_1(M)$ est un
groupe fini, alors toute hypersurface compacte $H$ est de type
fini.\end{cor} Preuve du corollaire 0-8 :

-Si $\pi_1(M)$ est un groupe fini, alors un th\'eor\` eme
d'Haefliger assure qu'il n'existe pas sur $M$ un feuilletage r\'eel
analytique (voir [9], chapitre 6).
\bigskip
-Si $\pi_2(M)\not =0$, alors nous avons l'alternative suivante (car
$M$ est de dimension 3, voir [9], chapitre 7, th\' eor\` eme de
Novikov) pour un feuilletage $F$ de $M$ : soit il admet une
composante de Reeb, ou bien toutes les feuilles sont compactes.
Quoiqu'il en soit, il poss\` ede donc  au moins une feuille compacte
; maintenant, int\' egrons la forme symplectique sur cette feuille -
$\omega$ et $J$ \'etant adapt\'ees - d'un c\^ot\'e, cette
int\'egrale est strictement positive et de l'autre elle est nulle \`
a l'aide de la formule de Stokes ($\omega$ est exacte sur $M$).

Pour le cas de $\C^{2}$, il est inutile de faire des hypoth\`
eses topologiques sur $M$ pour \'eviter le cas feuillet\'ee : ceci
repose sur un argument de A.Glutsyuk (voir [3] et la section [5]).
Nous avons donc le th\'eor\` eme suivant :
\begin{theo}\label{thm:main}Soit $(\C^{2},\omega_0)$ avec $\omega_0$ la structure
symplectique standard sur $\C^{2}$. Alors si $J$ est une
structure presque complexe adapt\'ee \` a $\omega_0$, toute
hypersurface compacte est de type fini pour $J$.
\end{theo}

Les th\'eor\`emes 0-7 et 0-9 reposent sur un r\'esultat clef qui
assure le prolongement des germes de disques $J$-holomorphes
contenus dans une hypersurface de $\R^{4}$ d\'efinie au voisinage
de $0$, en une courbe $J$-holomorphe ferm\'ee dans ce voisinage
(voir th\'eor\` eme 4-2 de la section 4).

Pour $\omega$ une structure symplectique donn\'ee, il existe
toutjours des structures presques complexes adapt\'ees. Par exemple,
pour $\C^{2}$ muni de la structure symplectique standard alors la
multiplication par $i$ est adapt\'ee, et nous retrouvons donc, dans
le cas de la dimension complexe $2$, un r\'esultat de
Diederich-Fornaess ([8]).

Pour finir, citons un exemple d'application du th\'eor\`eme 0-7 : le
cotangent d'une surface r\'eelle analytique que l'on munit de sa
structure symplectique naturelle, v\'erifie les hypoth\` eses du
th\'eor\` eme.

Mentionnons quelques questions ouvertes : pour le corollaire 0-5,
on peut imaginer appliquer un th\'eor\` eme
de compacit\'e aux $\gamma_n$ mais au pr\'ealable il faut les
prolonger, ce qui est possible d'apr\`es le th\'eor\` eme 4-2, mais
qui ne nous permet pas de faire l'\'economie du th\'eor\`eme de
Cartan-Kahler. De plus, il faut encore s'assurer que l'aire des
disques prolong\'es est uniform\'ement born\'e; dans le cas
int\'egrable, ceci est possible en utilisant quelques r\'esultats
\'el\'ementaires sur la multiplicit\'e des applications holomorphes
(au moins dans le cas d'une hypersurface alg\'ebrique, voir [7]), ce
qui n'a pas d'\'equivalent dans le cadre non int\'egrable.

N\'eanmoins, il nous semble int\'eressant de prouver s'il existe ou
non un proc\'ed\'e effectif de reconstruction de $\gamma$ \` a la
limite, \`a partir des $\gamma_n$, comme dans le cas int\'egrable
(voir [7]).
\medskip
Enfin, il nous para\^\i t naturel d'envisager ces probl\`emes en
dimension plus grande - dans un premier temps, sous l'hypoth\`ese
$M$ pseudoconvexe -, la difficult\'e suppl\'e-mentaire r\'eside dans
le fait suivant : on ne peut appliquer tel quel le r\'esultat de
Goldschmidt (th\'eor\`eme 3-2, section 3)car les projections ne sont
plus surjectives, et donc le syst\`eme n'est plus formellement
int\'egrable (voir d\'efinition 3-1).

\bigskip
\centerline{Rappels et pr\'eliminaires}
\bigskip

Soit $D$ le disque unit\'e de $\R^{2}$, $u:D\rightarrow \R^{4}$,
$(x,y)\rightarrow (u_1(x,y),\cdots,u_4(x,y))$, $J$ une structure
presque complexe sur $\R^{4}$, $J_0$ la structure complexe standard
sur $\R^{2}$, on veut r\'esoudre le syst\` eme diff\'erentiel, $E$,
au voisinage de z\'ero, suivant :
$$\begin{cases}& d\rho_{u(x,y)} od u_{(x,y)}=0\\
 & du_{(x,y)}J_{0} -J(u(x,y)) du_{(x,y)}=0.\end{cases}$$
Soit $J_{0}^{k}$ l'espace des k-jets en $0$ ; on associe \` a $E$ le
sous-espace $\varepsilon$ de $J^{1}_{0}$ d\'efini par : $\theta\in
\varepsilon$ si et seulement si $\theta$ v\'erifie\ E en $0$. On
consid\` ere \'egalement les prolongements de $E$, $E^k$ $\forall
k\geq 0$ (voir [4], P.137-140), obtenus par d\'erivation d'ordre $k$
de $E$ :
$$\begin{cases}& D^{\sigma}\big(d\rho_{u(x,y)} o d u_{(x,y)}\big)=0\\
& D^{\sigma} \big(du_{(x,y)}\ J_{0} -J(u(x,y))

du_{(x,y)}\big)=0.\end{cases}$$ o\` u $\sigma=(\sigma_1,\sigma_2)$
avec $\vert\sigma\vert\leq k$, et on leur associe $\varepsilon^k$,
le sous-espace de $J^{k+1}_0$  des jets qui v\'erifient $E^{k}$ en
$0$.

Nous rappelons ici des  r\'esultats  obtenus  dans [1], dont nous
aurons besoin par la suite pour la d\'emonstration du th\'eor\`eme
0-4 :
\begin{theo}\label{thm:main} Consid\`erons $L^{p,q}$ :
$(\R^{2n})^{p+q+1}\rightarrow \R$, la forme de Levi
d'ordre $p,q$(voir [1], pour la d\'efinition pr\'ecise).

Alors les deux propositions suivantes sont \'equivalentes:

(1) Il existe un disque $J$-holomorphe, $u$, r\'egulier tangent \`a
$M$ en $0$ \`a l'ordre $k+2$.

(2) Il existe, $X$, champ de vecteur $J$-tangent \`a $M$ (ie une
section r\'eguli\`ere de $T^{J}(M)$) v\'erifiant : $\forall\ (p,q),\
p+q\leq k-1$, $L^{p,q}\big(X(0),X.X(0),\cdots,X^{p+q+1}(0)\big)=0.$

De plus, on peut choisir $X$ tel que $X(0)={\partial u\over \partial
x}(0),\cdots,X^{k+1}(0)={\partial^{k+1} u\over \partial
x^{k+1}}(0)$.
\end{theo}
\begin{defn}On d\'efinit la d\'eriv\'ee
$(p,q)$-i\`eme en $0$
 d'un champ de vecteurs $X$ sur  $\R^{4}$,comme \'etant le vecteur $X$
 d\'eriv\'e $p$-fois dans sa propre direction et  $q$-fois dans la direction
 $JX$   $$   D^{p,q}_{X}X(0)=(JX)^{q}X^{p}\cdot X(0)$$
 \end{defn}
\begin{defn}La collection de toutes les
d\'eriv\'ees
 $(p,q)$-i\`eme en $0$,pour $0\leq p+q\leq k$, d'un champ de vecteurs $X$ sur
 $\R^{4}$ est appel\'ee son jet d'ordre $k$ et not\'e par
 $j^{k}_{0}(X)$. Ceci d\'efinie une application   $$  \begin{cases}
 j^{k}_{0}:
&
\Gamma(T\R^{4}) \rightarrow
  J^{k}_{0}(\R^{4})=(\R^{4})^{\frac{(k+1)(k+2)}{2}}\\
 & X \rightarrow (X(0),\dots, D^{p,q}_{X}X(0),\dots, D^{0,k}_{X}X(0))\end{cases}$$ et le
 jet $\xi\in J^{k}_{0}(\R^{4})$ est dit r\'ealisable
 s'il appartient \` a  $J^{k}_{0}(M):=j^{k}_{0}(\Gamma(T^{J}M)$.\end{defn}
\begin{prop}Soit $\xi\in
J^{k+1}_{0}(\R^{4})$. On note
 $[\xi]_k$ sa $J^{k}_{0}(\R^{4})$-composante et
 $$(\xi_{k+1,0},\cdots,\xi_{p,q},\cdots,\xi_{0,k+1})$$ sa partie homog\` ene
 d'ordre $k+1$. Alors
$$ \exists X\in\Gamma (T^{J}M),\ j_{0}^{k+1}(X)=\xi$$
si et seulement si
$$\begin{cases} &  \exists X_1\in\Gamma (T^{J}M),\
j_{0}^{k}(X_1)=[\xi]_k \\ & \forall\ (p,q),\ p+q=k+1,\
\xi_{p,q}-[j_{0}^{k+1}(X_1)]_{p,q}\in T_{0}^{J} M.\end{cases}$$
\end{prop}
Autrement dit, si $\xi$ est un $k+1$-jet r\'ealisable, les
composantes $N$ et $JN$ de sa partie homog\`ene de degr\'es $k+1$,
o\`u $N$ est le vecteur normal \`a $M$ en $0$, sont enti\`erement
d\'etermin\'ees par son jet d'ordre $k$; mais les composantes
$J$-tangentes (toujours de la partie homog\`ene d'ordre $k+1$) sont
totalement libres. \begin{cor}Soit $u$ un disque $J$-holomorphe
r\'egulier au voisinage de $0$. Alors $u$ est tangent \` a $M$ en
$0$ \`a l'ordre $k+2$(soit $\rho(u(z))=0 (z^{k+2})$)si et seulement
si son jet d'ordre $k+1$ est r\'ealisable.
\end{cor}

\section{Prolongement du jet d'un disque tangent \`a $M$
\`a l'ordre $k+2$}

\begin{prop}Soient $0$ un point de type infini et
$\varepsilon^{k+1}$, pour tout $k\geq 0$, le prolongement d'ordre
$k+1$ du syst\`eme diff\'erentiel $E$. Alors la projection $\pi_{k}
 : \epsilon^{k+1}\rightarrow
\epsilon^{k}$ est surjective pour tout $k\geq 0$.
\end{prop}
Preuve de la proposition 1-1 :
\medskip
Soit $u$ un disque $J$-holomorphe avec un ordre de contact $k+2$
avec $M$; montrer que $\pi_k$ est surjective revient \`a prouver
l'existence de $v$, un disque $J$-holomorphe qui lui tangente \`a
l'ordre $k+3$, avec $J_{0}^{k+1}(u)=J_{0}^{k+1}(v)$.
\medskip
On sait qu'il existe $\omega$, et donc $J_{0}^{k+2}(\omega)$, tel
que $\omega$ \`a un ordre de contact $k+3$ avec $M$ ; d'apr\`es le
th\'eor\`eme 0-10, on a par cons\'equent,
$$L^{p,q}\big(\omega(0),{\partial \omega\over
\partial x}(0),\cdots,{\partial^{p+q+1}\omega\over \partial
x^{p+q+1}}(0)\big)=0,$$ $\forall p+q\leq k$. Si $\phi$ est une
fonction holomorphe au voisinage de $0$ avec $\phi(0)=0$, alors
$\omega o\phi$ a \'egalement un ordre de contact $k+3$ avec $M$, et
donc pour tout $Y$ un champ de vecteurs $J$-tangent \`a $M$ de la
forme $(\alpha+\beta J)X_{\omega}$(o\`u $X_{\omega}$ est un champ de
vecteurs associ\'e \`a $\omega$ par le th\'eor\`eme 0-10, et
$\alpha$, $\beta$ sont des fonctions r\'eelles r\'eguli\`eres ),
nous avons : $\forall (p,q)$ avec $p+q\leq k$, $L^{p,q}(Y)(0)=0.$
\medskip
Remarquons que $dim_{\C}(T^J)=1$, et donc
$X_u=(\tilde{\alpha}+\tilde{\beta}J)X_{\omega}$, ce qui entra\^\i
ne, d'apr\`es ce qui pr\'ec\`ede, $\forall (p,q)$ avec $p+q\leq k$,
$L^{p,q}(X_{u})(0)=0.$ A nouveau, le th\'eor\`eme 0-10 prouve
l'existence d'un disque, $v$, $J$-holomorphe ayant un ordre de
contact $k+3$ avec $M$ en $0$; de plus,
$J_{0}^{k+1}(v)=J_{0}^{k+1}(u)$ et ${\partial^{k+2}v\over \partial
x^{k+2}}(0)=X_{u}^{k+2}(0)$.

\section{Symbole
d'un syst\`eme diff\'erentiel}

Pour une description plus d\'etaill\'e, on pourra consulter [4] et
[5]. Afin d'\'eviter des notations un peu lourdes, nous nous
contentons ici de rappeler uniquement ce dont nous aurons besoin.
Consid\'erons un syst\`eme d'ordre $1$ d\'efini par :
$$F_l(x,y,u_i(x,y),\partial_1 u_i,\partial_2 u_i)=0, \forall l\in
\{1\cdots,m\}.$$ Pour $k\geq 1$, on d\'efinit le symbole de
$\varepsilon^k$, $g^k$, qui \`a $\theta\in \varepsilon^k$ associe
l'espace vectoriel
$$g^k(\theta):=\big\{(\zeta_{k+1,0},\cdots,\zeta_{p,q},\cdots,\zeta_{0,k+1})
/\sum_{i}\zeta^{i}_{p+1,q-1}\partial^{i}_{1}F_l(\theta)\zeta^{i}_{p-1,q+1}
\partial^{i}_{2}F_l(\theta)=0\big\}.$$ Ici, nous devons
pr\'eciser que nous consid\'erons $F$ comme une fonction des
variables ind\'ependantes $(x,y,u_i,p^{i}_1,p^{i}_2),$ et
$\partial^{i}_{1}F_l$, $\partial^{i}_{2}F_l$ signifient
respectivement ${\partial F_l\over \partial p^{i}_1}$, ${\partial
F_l\over \partial p^{i}_2}$. Il est clair que $g^k(\theta)$ est
enti\`erement d\'etermin\'e par la partie homog\`ene d'ordre $1$ de
$\theta$; pour cette raison, on consid\`ere $\theta\in
\varepsilon\rightarrow g^k(\theta)$.
\medskip
Dans notre cas, nous savons (voir section pr\'ec\`edente) que les
projections $\pi_{k} :\varepsilon^{k+1}\rightarrow\varepsilon^k$
sont surjectives, ce qui entra\^\i ne que $g^k(\theta)$ (pour
$\theta\in \varepsilon$) peut \^etre d\'ecrit de la mani\`ere
suivante :
$$g^k(\theta):=\big\{\big((T_1-T_2)_{k+1,0},\cdots,(T_1-T_2)_{p,q},\cdots,(T_1-T_2)_{0,k+1}\big)\big\}$$
$$\hbox{o\`u}\ (T_1,T_2)\in \varepsilon^k\times\varepsilon^k,\
J^{k}_0(T_1)=J^{k}_0(T_2)\ \hbox{et}\
J^{1}_0(T_1)=J^{1}_0(T_2)=\theta$$ et
$$\big((T_1-T_2)_{k+1,0},\cdots,(T_1-T_2)_{p,q},\cdots,(T_1-T_2)_{0,k+1}\big)$$
d\'esigne la partie homog\` ene d'ordre $k+1$ de $T_1-T_2$.
D'apr\`es le corollaire 0-14, ceci signifie que $T_1$ et $T_2$ sont
des jets d'ordre $k+1$ r\'ealisables ayant la m\^eme composante
d'ordre $k$, nous pouvons appliquer alors la proposition 0-13 : la
partie homog\`ene de $T_1-T_2\in (T^{j}M)^{k+1}$. Nous venons donc
de montrer $(\theta,g^k(\theta))$ est un fibr\'e vectoriel de fibre
$(T^JM)^{k+1}$ au dessus de $\varepsilon$.

\section{Preuve du th\'eor\`eme 0-4}

Nous avons $\R^{4}$ muni d'une structure presque complexe $J$
d\'ependant analytiquement du point base et $M$ une hypersurface
r\'eelle analytique, par cons\'equent,le syst\`eme $E$ de la partie
1 est \`a coefficients analytiques ; de plus d'apr\`es [5], $E$
v\'erifie: \begin{defn}On dit qu'un syst\`eme
diff\'erentiel $E$ est formellement int\'egrable si et seulement si
$(\theta,g^{k+1}(\theta))$ est un fibr\'e vectoriel pour tout $k\geq
0$ au dessus de $\varepsilon$ et si de plus l'op\'erateur de
restriction $\pi_{k}:\varepsilon^{k+1}\rightarrow \varepsilon^{k}$
est surjectif.
\end{defn}
\begin{theo}\label{thm:main}Soit $E$ un syst\`eme diff\'erentiel
analytique formellement int\'egrable, alors pour tout $\theta \in
\varepsilon^{k}$, il existe une solution analytique, $u$ de $E$,
telle que $J_{0}^{k+1}(u)=\theta$.
\end{theo}

Preuve du th\'eor\` eme 0-4 :
\medskip
le point $0$ est un point de type infini ce qui entra\^\i ne la
formelle int\'egrabilit\'e du syst\`eme $E$, d'apr\`es les parties
$1$ et $2$ ; le th\'eor\` eme 0-4 est alors la cons\'equence directe
du th\'eor\`eme 3-2. \begin{rem}Il est clair que m\^eme si
l'hypersurface et la structure sont seulement de classe
$C^{\infty}$, le syst\`eme $E$ reste formellement
int\'egrable.\end{rem}

\section{Preuve du th\'eor\`eme 0-7}

Dans cette partie, on convient de noter par une majuscule l'image
g\'eom\'etrique d'un disque  et par une minuscule le disque
proprement dit.

Nous allons commencer par d\'efinir le prolongement d'un disque
$J$-holomorphe: \begin{defn}Soit $u:D\rightarrow
\R^{4}$, $0\rightarrow 0$ un germe de disque $J$-holomorphe dont
l'image est incluse dans $\Omega$, un voisinage ouvert de $0$ dans
$\R^{4}$; on dit que $u$ admet un prolongement \`a $\Omega$ ou se
prolonge \` a $\Omega$, s'il existe $w$ une courbe $J$-holomorphe
ferm\'ee dans $\Omega$ dont l'image $W$ contient $U$, l'image de
$u$.\end{defn} Commencons par d\'emontrer le th\'eor\`eme suivant:
\begin{theo}\label{thm:main}Soit $0\in M$ une hypersurface
r\'eelle analytique de $\R^{4}$ d\'efinie sur $\Omega$; soit $u$
un germe de disque $J$-holomorphe,r\'egulier, passant par $0$ et
contenu dans $M$. Alors, si $M$ n'est pas feuillet\'ee en courbe
$J$-holomorphe, $u$ se prolonge \` a $\Omega$ en une courbe
r\'eguli\`ere $J$-holomorphe.
\end{theo}

Preuve du th\'eor\`eme 4-2:

Soit $x_0\in \overline{U}$, on va montrer que $u$ se prolonge en une
courbe $J$-holomorphe au voisinage de $x_0$. Consid\'erons, $X$ un
champ de vecteurs $J$-tangent \` a $M$ au voisinage de $x_0$ et
$L_{\zeta}(X(\zeta),X(\zeta)):=L_{\zeta}$ la forme de levi en
$\zeta\in M$ appliqu\'ee au vecteur $J$-tangent $X(\zeta)$.
Clairement $U\subset\{L_{\zeta}(X(\zeta),X(\zeta))=0\}$; $M$
n'\'etant pas feuillet\'ee, $\{L_{\zeta}(X(\zeta),X(\zeta))=0\}$ est
de dimension $2$, et  $U$ dans l'une des composantes irr\'eductibles
de cet ensemble r\'eel
   analytique, disons $\{g(\zeta)=0\}$ ($g$ \'etant l'un des facteurs
   irr\'eductibles en $x_0$ de $L_{\zeta}$);d'apr\`es le corollaire 0-5, il existe
   $\tilde{u}$, un germe de disque $J$-holomorphe en $x_0$, dont l'image $\tilde{U}$ est
contenue dans $M$ ; nous allons montrer qu'en fait $\tilde
{U}\subset  \{g(\zeta)=0\}$. Nous avons
   vu dans la section 2, que si $\theta$ est un germe de disque $J$-holomorphe
   en $Y$ qui a un ordre de contact $k+2$ avec $M$, alors on  peut trouver un
   re-param\'etrage, $\alpha^{\theta}_{k}$ tel que
   ${\partial^{i}(\theta(\alpha^{\theta}_{k}))\over \partial
   x^{i}}(0)=X^{i}(Y)$, pour tout $i\leq k$. Soit $k\in\N$ quelconque et
   $(x_n)$ une suite de points de $U$ qui tend vers $x_0$, on note $u_{x_n}$ le
   germe de disque en $x_n$ d\'efini par $u$; choisissons
   $\alpha^{u_{x_n}}_{k}$ un re-param\'etrage ayant la propri\'et\'e
   ci-dessus; d'apr\`es ce qui pr\'ec\` ede, en passant \`a la limite sur les
   jets d'ordre $k$ des disques $u_{x_n}(\alpha^{u_{x_n}}_{k})$, on obtient
   que le jet d'ordre $k$ compatible, dont les d\'eriv\'ees d'ordre inf\'erieur
   \`a $k$ en $x$ sont donn\'ees par $X^{i}(x_0)$, annule \`a l'ordre $k$ la
   fonction $g$. D'autre part, toujours en accord avec ce qui pr\'ec\` ede, je peux trouver un re-param\'etrage
   $\alpha^{\tilde{u}}_{k}$ tel que${\partial^{i}(\tilde{u}(\alpha^{\tilde{u}}_{k}))\over \partial
   x^{i}}(0)=X^{i}(x_0)$, pour tout $i\leq k$, ce qui entra\^\i ne
   $g(\tilde{u}(\alpha^{\theta}_{k}) )$ s'annule \`a l'ordre $k$ en $x_0$, et par
   cons\'equent $g(\tilde{u})$ \'egalement. Je peux effectuer cette op\'eration
   pour tout $k$, ce qui implique $\tilde{U}\subset \{g=0\}$ et $u$ se prolonge
   en une courbe r\'eguli\`ere au voisinage de $x_0$.

 En prolongeant de proche en proche de la mani\`ere pr\'ec\` edente,on
obtient le r\'esultat en s'assurant n\'eanmoins que si l'on revient
en $x_0$, le germe de disque obtenu reste dans la m\^eme composante
irr\'eductible de $L_{\zeta}=0$, soit $\{g(\zeta)=0\}$ ; ceci est
assur\'e car le jet du disque en question (quitte \` a
re-param\'etrer) est encore donn\'e par le jet d'ordre $k$
compatible, dont les d\'eriv\'ees d'ordre inf\'erieur
   \`a $k$ en $x$ sont les $X^{i}(x_0)$, qui annule $g$, d'apr\` es ce qui
pr\'ec\` ede.

A pr\'esent, nous pouvons justifier l'affirmation suivante : (voir
fin de la preuve du corollaire 0-5) \begin{rem}$M$ ne
peut contenir un germe de disque singulier en $0$.
\end{rem}
Si $L_{\zeta}$ est identiquement nulle, alors $M$ est feuillet\'ee
en courbe r\'eguli\` ere et l'affirmation est \'etablie.
\medskip
Dans le cas contraire, $\{L_{\zeta}=0\}$ est de dimension $2$ en
$0$; raisonnons par l'absurde, supposons qu'il existe $u$ singulier
en $0$ dont l'image $U\subset M$; alors $U\subset \{g=0\}$ o\`u
$\{g=0\}$ est la composante irr\'eductible de $L_{\zeta}$ qui
contient $U$. Les points singuliers d'un disque $J$-holomorphe sont
isol\'es, ce qui entra\^\i ne qu'il existe une suite $(a_n)$ qui
tend vers $0$ avec $u$ qui est un germe de disque r\'egulier en
$a_n$; la preuve du th\'eor\`eme 4-2 prouve alors que $\{g=0\}$
contient un germe de disque r\'egulier en $0$ et donc, est une
vari\'et\'e  r\'eelle analytique de dimension $2$ au voisinage de
$0$, ce qui implique que $u$ est r\'egulier en $0$. Nous avons
\'egalement une version globale de cette affirmation :
\begin{rem}$M$ ne peut contenir une courbe
$J$-holomorphe singuli\` ere (attention, $M$ pourrait contenir des
courbes singuli\` eres sans contenir de germes de disques singuliers
: il suffit que les composantes irr\'eductibles de la courbe en un
point multiple soient des disques r\'eguliers).
\end{rem}
En effet, supposons qu'il existe $C$ une courbe singuli\` ere dans
$M$, alors $C$ contient un germe de disque r\'egulier qui est
lui-m\^eme contenu dans une courbe r\'eguli\` ere, $C^{'}$,d'apr\`es
le th\'eor\` eme 4-2; mais alors $C$ et $C^{'}$ sont identiques, et
donc $C$ n'a pas de singularit\'es.
\medskip
{Remarques : }

Supposons que $M$ soit une hypersurface r\'eelle analytique
plong\'ee(pas forc\'ement compacte)dans une vari\'et\'e de dimension
4. Alors si $M$ est non feuillet\'ee, la preuve du th\'eor\` eme 4-2
entra\^\i ne que tout germe de disque  inclu dans $M$ se prolonge en
une courbe ferm\'ee dans $M$. Maintenant si $M$ est feuillet\'ee et
simplement connexe, alors le germe de disque se prolonge par l'une
des feuilles, mais toutes les feuilles d'un feuilletage
transversalement r\'eel analytique sont ferm\'ees([9], chapitre 6),
et donc le germe se prolonge en une courbe ferm\'ee \'egalement. On
a prouv\'e la proposition suivante : \begin{prop} Soit
$M$ une hypersurface d' une vari\'et\'e presque complexe $(V,J)$. Si
$M$ est simplement connexe, alors tout germe de disque contenu dans
$M$ se prolonge en une courbe ferm\'ee dans $M$
\end{prop}

\bigskip
 \centerline{Fin de la preuve du th\'eor\` eme 0-7}
\bigskip

Si $M$ n'est pas feuillet\'ee, la preuve du th\'eor\` eme 4-2 assure
que $U\subset C$ avec $c$ une courbe ferm\'ee r\'eguli\`ere dans
$M$, ce qui entra\^\i ne que  $C$ est compacte sans bord. D'un
c\^ot\'e - $J$ et $\omega$ \'etant adapt\'ees - $\int_{C}\omega\geq
0$ et de l'autre  $\int_{C}\omega =\int_{C}d\xi =0$ en utilisant la
formule de Stokes; ce qui fournit une contradiction.

\section{Preuve du th\'eor\` eme 0-9}

Nous allons reproduire l'argument de A.Glutsyuk (voir aussi [3]).
D'apr\` es le Th\'eor\` eme 0-7, il suffit de montrer que $M$ n'est
pas feuillet\'ee. Supposons que ce soit le cas :

le feuilletage ne d\'epend que de la structure $J$ sur une boule
contenant $M$ car $M$ est compacte. Construisons $\tilde J$ une
structure presque complexe sur $\C P^{2}$ adapt\'ee \` a $\omega_0$
telle que: $\tilde J$ est la structure donn\'ee sur cette boule et
$\tilde J=J_0$, la structure complexe habituelle sur l'hyperplan \`
a l'infini (la nouvelle structure n'est peut \^etre pas  r\'eelle
analytique  mais on peut le faire dans la cat\'egorie $C^\infty$, ce
qui suffit pour ce qui suit).

Soit $\Omega$, l'ensemble des courbes rationnelles $\tilde
J$-holomorphes rencontrant $T^{\tilde J}M$ et   homologues \` a une
droite projective $\C P^{1}\subset\C P^{2}$(les courbes $\tilde
J$-holomorphes homologues \`a $\C P^{1}\subset \C P^{2}$ sont
appel\'ees les  $\tilde J$-droites); montrons que $\Omega$ est
ouvert et ferm\'e dans l'espace des $\tilde J$-droites.

1)$\Omega$ ferm\'e: soit $D_N$ une suite de droites rencontrant
$T^{\tilde J}M$  avec $D_N$ qui tend vers $D$, et $a_N$ une suite de
points de $M$ tels que $D_N$ rencontre  $T^{\tilde J}M$ en $a_N$;
  si $a$ une  valeur d'adh\'erence de la suite $a_N$ ($M$ est compacte),
  alors $a\in D$ et donc $D$ rencontre  $T^{\tilde J}M$ en $a$.

2)$\Omega$ est ouvert : soit $D_0$ une droite de $\Omega$ et $a$ un
point de $M$
  tel que $D_0$ rencontre  $T^{\tilde J}M$ en $a$. $F_a$, la feuille passant par
  $a$, est ferm\'ee dans une boule suffisamment petite autour de
  $a$. Si $D$ est  une petite perturbation de $D_0$ alors l'indice
d'intersection en  $a$, de $F_a$ et $D$, est le m\^eme que celui de
$D_0$ et $F_a$, donc strictement positif, ce qui entra\^\i ne : $D$
rencontre le feuilletage et  $\Omega$ est ouvert.

D'apr\`es un r\'esultat de Gromov ([6]), par deux points distincts
de $\C P^{2}$, il passe une unique $\tilde J$-droite; de plus, si
l'on note $D(a,b)$ la droite passant par $(a,b)\in \C P^{2}\times\C
P^{2}$, l'application $(a,b)\rightarrow D(a,b)$ est r\'eguli\` ere
et donc l'ensemble des  $\tilde J$-droites est connexe. Par
cons\'equent, $\Omega$ est l' ensemble des $\tilde J$-droites , ce
qui produit une contradiction car la $\tilde J$ droite \` a l'infini
ne rencontre pas  $T^{\tilde{J}}M$, $M$ \'etant compacte dans
$\C^{2}$ et $\tilde J=J_0$ \` a l'infini.

\end{document}